\documentclass[12pt,reqno]{amsart}
\usepackage{amssymb,latexsym,amsmath,mathrsfs}
\usepackage{hyperref}  
\usepackage[matrix, arrow]{xy}
\xyoption{all} 

\newcommand{\K}{\mathcal{K}}

\newcommand{\U}{\mathcal{U}}
\newcommand{\C}{\mathbb{C}}
\newcommand{\N}{\mathbb{N}}
\newcommand{\T}{\mathbb{T}}
\newcommand{\Z}{\mathbb{Z}}

\newcommand{\Oinf}{{\mathcal{O}_\infty}}
\newcommand{\Td}{{\mathcal{T}_2}}
\newcommand{\Tt}{{\mathcal{T}_3}}
\newcommand{\Tn}{{\mathcal{T}_n}}

\newcommand{\pf}{\noindent {\mbox{\textit{Proof}. }} }
\newcommand{\ie}{\textit{i.e.\/}\ } 

\newcommand{\cst}{C$^*$} 

\newtheorem{thm}{Theorem}[section]
\newtheorem{cor}[thm]{Corollary}
\newtheorem{lem}[thm]{Lemma}
\newtheorem{prop}[thm]{Proposition}

\theoremstyle{definition}
\newtheorem{defi}[thm]{Definition}

\newtheorem{rem}[thm]{Remark}
\newtheorem{rems}[thm]{Remarks}

\theoremstyle{remark}

\numberwithin{equation}{section}
\title{$K_1$-injectivity for properly infinite \cst-algebras}

\author{\'Etienne Blanchard}

\subjclass[2010]{Primary: 46L80; Secondary: 46L06, 46L35} 

\keywords{$K_1$-injectivity, Proper Infiniteness, \cst-algebras}

\begin{document} 
\maketitle 
\centerline{\textit{Dedicated to Alain Connes on the occasion of his 60th birthday.}}
\section{Introduction} 
One of the main tools to classify \cst-algebras is the study of its projections and its unitaries. 
It was proved by Cuntz in \cite{Cu81} that if $A$ is a \textit{purely infinite} simple \cst-algebra, 
then the kernel of the natural map for the unitary group $\U(A)$ 
to the $K$-theory group $K_1(A)$ 
is reduced to the connected component $\U^0(A)$, \ie $A$ is \textit{$K_1$-injective} (see~\S 3). 
We study in this note a finitely generated \cst-algebra, the $K_1$-injectivity of which would imply the $K_1$-injectivity of all unital \textit{properly infinite} \cst-algebras. \\ \indent 
Note that such a question was already considered in \cite{Blac}, \cite{BRR}. 

\medskip 
The author would like to thank H. Lin, R. Nest, M. R\o rdam and W. Winter for helpful comments.

\section{Preliminaries} 

Let us first review briefly 
the theory introduced by Cuntz (\cite{Cu78}) of comparison of positive elements in a \cst-algebra. 
\begin{defi}(\cite{Cu78}, \cite{Ro92}) 
Given two positive elements $a, b$ in a \cst-algebra $A$, one says that: \\
-- $a$ is \textit{dominated} by $b$ (written $a \precsim b$) if and only if there is a sequence 
$\{ d_k; k\in\N\}$ in~$A$ such that $\| d_k^*b d_k^{} - a \| \to 0$ when $k\to\infty$, \\
-- $a$ is \textit{properly infinite} if $a\not= 0$ and 
$a\oplus a \precsim a \oplus 0$ in the \cst-algebra $M_2(A):=M_2(\C)\otimes A$. 
\end{defi}

This leads to the following notions of infiniteness for \cst-algebras. 
\begin{defi} (\cite{Cu78}, \cite{Cu81}, \cite{KiRo00}) 
A unital \cst-algebra $A$ is said to be:\\
-- \textit{properly infinite} if its unit $1_A$ is properly infinite in $A$, \\
-- \textit{purely infinite} if all the non zero positive elements in $A$ are properly infinite in~$A$.
\end{defi} 
\begin{rem} Kirchberg and R\o rdam have proved in \cite[Theorem 4.16]{KiRo00} that 
a \cst-algebra $A$ is purely infinite (in the above sense) if and only if 
there is no character on the \cst-algebra $A$ and 
any positive element $a$ in $A$ which lies in the closed two-sided ideal generated by another positive element $b$ in $A$ satisfies $a\precsim b$. 
\end{rem}

The first examples of such \cst-algebras were given by Cuntz in \cite{Cu81}: 
For any integer $n\geq 2$, the \cst-algebra $\Tn$ is the universal unital \cst-algebra generated by $n$ isometries $s_1, \ldots, s_n$ satisfying the relation 
\begin{equation}\label{Tn} 
s_1^{}s_1^* +\ldots+s_n^{}s_n^*\leq 1
\end{equation} 
Then, the closed two sided ideal in $\Tn$ generated by the \textit{minimal} projection 
$p_{n+1}:=1-s_1^{}s_1^*-\ldots -s_n^{}s_n^*$ is isomorphic to the \cst-algebra $\K$ of compact operators on an infinite dimension separable Hilbert space and one has an exact sequence 
\begin{equation}\label{TnOn} 
0\to \K \to \Tn\mathop{\longrightarrow}\limits^\pi \mathcal{O}_n\to 0\,,
\end{equation}
where the quotient $\mathcal{O}_n$ is a purely infinite \textit{simple} unital nuclear \cst-algebra (\cite{Cu81}). 

\begin{rem}\label{rem2.4} 
A unital \cst-algebra $A$ is properly infinite if and only if 
there exists a unital $*$-homomorphism $\Td\to A$. 
\end{rem}


\section{$K_1$-injectivity of $\Tn$}
Given a unital \cst-algebra $A$ with unitary group $\U(A)$, 
denote by $\U^0(A)$ the connected component of $1_A$ in $\U(A)$. 
For each strictly positive integer $k\geq 1$, the upper diagonal embedding 
$u\in\U(M_k(A)\,)\mapsto (u\oplus 1_A)\in\U(M_{k+1}(A)\,)$ 
sends the connected component $\U^0(M_k(A)\,)$ into $\U^0(M_{k+1}(A)\,)$, 
whence a canonical homomorphism $\Theta_A$ from $\U(A)\big{/}\U^0(A)$ to $K_1(A):=\mathop{\lim}\limits_{k\to\infty}\, \U(M_k(A)\,)\big{/} \U^0(M_k(A)\,)$. 
As noticed by Blackadar in \cite{Blac}, this map is (1) neither injective, (2) nor surjective in general: 
\begin{enumerate} 
\item If $\mathfrak{U}_2$ denotes the compact unitary group of the matrix \cst-algebra $M_2(\C)$, 
$A:=C(\mathfrak{U}_2\times \mathfrak{U}_2, M_2(\C)\,)$ and $u, v\in \U(A)$ are the two unitaries
$u(x, y)=x$ and $v(x, y)=y$, then $z:=uvu^*v^*$ is not unitarily homotopic to $1_A$ in $\U(A)$ 
but the unitary $z\oplus 1_A$ belongs to $\U^0(M_2(A) )$ (\cite{AJT}). 
\item If $A=C(\mathbb{T}^3)$, then $\U(A)/\U^0(A)\cong\Z^3$ but $K_1(A)\cong\Z^4$. 
\end{enumerate}

\begin{defi} 
The unital \cst-algebra $A$ is said to be \textit{$K_1$-injective} if the map $\Theta_A$ is injective. 
\end{defi} 

Cuntz proved in \cite{Cu81} that $\Theta_A$ is surjective as soon as the \cst-algebra $A$ is properly infinite and 
that it is also injective if the \cst-algebra $A$ is simple and purely infinite. 
Now, the $K$-theoretical six-term cyclic exact sequence associated to the exact sequence (\ref{TnOn}) implies that $K_1(\Tn)=0$ since $K_1(\K)=K_1(\mathcal{O}_n)=0$. 
Thus, the map $\Theta_{\Tn}$ is zero. 
\begin{prop}\label{prop3.2} 
For all $n\geq 2$, the \cst-algebra $\Tn$ is $K_1$-injective, 
\ie any unitary $u\in\U(\Tn)$ is unitarily homotopic to $1_{\Tn}$ 
in $\U(\Tn)$ (written $u\sim_h 1_{\Tn}$). 
\end{prop}
\pf 
The \cst-algebras $\Tn$ have real rank zero by Proposition 2.3 of \cite{Zha}. 
And Lin proved that 
any unital \cst-algebra of real rank zero is $K_1$-injective (\cite[Corollary~4.2.10]{Lin}). 
\qed

\begin{cor}\label{cor3.3} 
If $\alpha:\Tt\to\Tt$ is a unital $*$-endomorphism, 
then its restriction to the unital copy of $\Td$ generated by the two isometries $s_1, s_2$ 
is unitarily homotopic to $id_\Td$ among all unital $*$-homomorphisms $\Td\to\Tt$. 
\end{cor} 
\pf The isometry $\sum_{k=1, 2} \alpha(s_k) s_k^*$ extends 
to a unitary $u\in\U(\Tt)$ such that $\alpha(s_k)=us_k$ for $k=1, 2$ (\cite[Lemma 2.4]{BRR}). 
But Proposition \ref{prop3.2} yields that $\U(\Tt)=\U^0(\Tt)$, 
whence a homotopy $u\sim_h 1$ in $\U(\Tt)$, and so the desired result holds. 
\qed 

\begin{rem}\label{rem3.4} 
The unital map $\iota: \mathbb{C}\to\Td$ induces an isomorphism in $K$-theory. 
Indeed, $[1_{\Td}]=[s_1s_1^*] + [s_2s_2^*] +[p_3]=2\, [1_{\Td}] + [p_3]$ 
and so $[1_{\Td}]=-[p_3]$ is invertible in $K_0(\Td)$. 
\end{rem} 

\section{$K_1$-injectivity of properly infinite \cst-algebras} 

Denote by $ \Td\ast_\C\Td$ the universal unital free product with amalgamation over $\C$ 
(in the sequel called full unital free product) 
of two copies of $\Td$ amalgamated over $\C$ 
and let $\jmath_0$, $\jmath_1$ be the two canonical unital inclusions of $\Td$ in $ \Td\ast_\C\Td$. 
We show in this section that the $K_1$-injectivity of $ \Td\ast_\C\Td$ is equivalent to 
the $K_1$-injectivity of all the unital properly infinite \cst-algebras. 
The proof is similar to that 
of the universality of the full unital free product $\Oinf\ast_\C\Oinf$ 
(see Theorem 5.5 of \cite{BRR}). 

\begin{defi} (\cite{Blan}, \cite[\S 2]{BRR}) 
If $X$ is a compact Hausdorff space, 
a unital $C(X)$-algebra is a unital \cst-algebra $A$ endowed with a unital $*$-homomorphism from the \cst-algebra $C(X)$ of continuous functions on $X$ to the centre of $A$. 

For any non-empty closed subset $Y$ of $X$, we denote by $\pi_Y^A$ 
(or simply by $\pi_Y$ if no confusion is possible) 
the quotient map from $A$ to the quotient $A_Y$ 
of $A$ by the (closed) ideal $C_0(X\setminus Y) \cdot A\,$. 
For any point $x\in X$, we also denote by $A_x$ the quotient $A_{\{ x\}}$ and by $\pi_x$ the quotient map $\pi_{\{ x\}}$. 

\end{defi}

\begin{prop}\label{main} 
The following assertions are equivalent. 
\begin{enumerate}
\item[(i)] $ \Td\ast_\C\Td$ is $K_1$-injective. 
\item[(ii)] $\mathcal{D}\!:=\!\!\{f\in C([0, 1] ,  \Td\ast_\C\Td)\,;\, f(0)\in\jmath_0(\Td)\,\mathrm{and}\,f(1)\!\in\!\jmath_1(\Td)\,\}$ is properly infinite. 
\item[(iii)] There exists a unital $*$-homomorphism $\theta: \Td\to\mathcal{D}$. 
\item[(iv)] There exists a projection $q\in\mathcal{D}$ with 
$\pi_0(q)=\jmath_0(s_1s_1^*)$ and $\pi_1(q)=\jmath_1(s_1s_1^*)\,$. 
\item[(v)] Any unital properly infinite \cst-algebra $A$ is $K_1$-injective. 
\end{enumerate}
\end{prop} 
\pf (i)$\Rightarrow$(ii) We have a pull-back diagram
\begin{equation*} 
\xymatrix{& \mathcal{D} \ar[ld] \ar[rd] & \\ \mathcal{D}_{[0,\frac{1}{2}]} \ar[rd]_{\pi_{\frac{1}{2}}}
 && \mathcal{D}_{[\frac{1}{2},1]} \ar[ld]^{\pi_{\frac{1}{2}}} \\ &  \Td\ast_\C\Td
 &}
\end{equation*} 
and the two \cst-algebras $\mathcal{D}_{[0,\frac{1}{2}]}$ and $\mathcal{D}_{[\frac{1}{2},1]}$ are properly infinite (Remark \ref{rem2.4}). 
Hence, the implication follows from \cite[Proposition 2.7]{BRR}. 

\medskip 
(ii)$\Rightarrow$(iii) is Remark~\ref{rem2.4} applied to the \cst-algebra $\mathcal{D}$. 

\medskip 
(iii)$\Rightarrow$(iv) The two full, properly infinite projections $\jmath_0(s_1s_1^*)$ and $\pi_0\circ\theta(s_1s_1^*)$ are unitarily equivalent in $\jmath_0(\Td)$ by \cite[Lemma 2.2.2]{LLR} 
and \cite[Proposition~2.3]{BRR}. 
Thus, they are homotopic among the projections in the \cst-algebra $\jmath_0(\Td)$ 
(written $\jmath_0(s_1s_1^*)\sim_h\pi_0\circ\theta(s_1s_1^*)\,$) by Proposition~\ref{prop3.2}. 
Similarly, $\pi_1\circ\theta(s_1s_1^*)\sim_h \jmath_1(s_1s_1^*)$ in $\jmath_1(\Td)$. 
Further, $\pi_0\circ\theta(s_1s_1^*)\sim_h \pi_1\circ\theta(s_1s_1^*)$ in $ \Td\ast_\C\Td$ by hypothesis, 
whence the result by composition. 

\medskip 
(iv)$\Rightarrow$(v) By \cite[Proposition 5.1]{BRR}, it is equivalent to prove that 
if $p$ and $p'$ are two properly infinite full projections in $A$, 
then there exist full properly infinite projections $p_0,$ and $p_0'$ in $A$ such that 
$p_0\leq p$, $p_0'\leq p'$ and $p_0\sim_h p_0'\,$. 

Fix two such projections $p$ and $p'$ in $A$. 
Then, there exist unital $*$-homomorphisms 
$\sigma: \Td\to pAp$, $\sigma': \Td\to p'Ap'$ and 
isometries $t, t'\in A$ such that $1_A =t^* p t=(t')^* p' t'\,$. 
Now, one thoroughly defines unital $*$-homomorphisms $\alpha_0:\Td\to A$ and $\alpha_1:\Td\to A$ by 
\begin{center} 
$\alpha_0(s_k):= \sigma(s_k) \cdot t$\quad and\quad 
$\alpha_1(s_k):= \sigma'(s_k) \cdot t'$ \qquad for $k=1, 2\,$, 
\end{center} 
whence a unital $*$-homomorphism $\alpha:=\alpha_0\ast\alpha_1:  \Td\ast_\C\Td\to A$ 
such that $\alpha\circ\jmath_0=\alpha_0$ and $\alpha\circ\jmath_1=\alpha_1\,$. 
\\ \indent 
The two full properly infinite projections $p_0=\alpha_0(s_1s_1^*)$ and $p_0'=\alpha_1(s_1s_1^*)$ satisfy $p_0\leq p$ and $p_0'\leq p'\,$. 
Further, the projection $(id\otimes\alpha)(q)$ gives a continuous path of projections in $A$ from $p_0$ to $p_0'\,$. 
\qed
\begin{rem}
The \cst-algebra $M_2(\mathcal{D})$ is properly infinite by \cite[Proposition~2.7]{BRR}. 
\end{rem} 

\begin{lem}\label{K-Td}  
$K_0( \Td\ast_\C\Td)=\Z$ and $K_1( \Td\ast_\C\Td)=0$
\end{lem}
\pf The commutative diagram 
$\raisebox{.8cm}{
\xymatrix{\C \ar[d]_{\imath_0} \ar[r]^{\imath_1} & \Td \ar[d]^{\jmath_1}\\
 \Td \ar[r]^{\jmath_0\quad} &  \Td\ast_\C\Td} }
$
yields by \cite[Theorem 2.2]{G} a six-term cyclic exact sequence 
$$\begin{array}{ccccc}
K_0(\C)=\Z&\mathop{\longrightarrow}\limits^{(\imath_0\oplus \imath_1)_*} & 
K_0(\Td\oplus\Td)=\
\Z\oplus\Z&\mathop{\longrightarrow}\limits^{(\jmath_0)_*-(\jmath_1)_*}&
K_0( \Td\ast_\C\Td)\\
\uparrow&&&&\downarrow\\
K_1( \Td\ast_\C\Td)&\longleftarrow &K_1(\Td\oplus\Td)=0\oplus 0&\longleftarrow&K_1(\C)=0
\end{array}$$ 
Now, Remark~\ref{rem3.4} implies that the map $(\imath_0\oplus \imath_1)_*$ is injective, whence the equalities. \qed 
\begin{rem} 
Skandalis noticed that the \cst-algebra $\Td$ is $KK$-equivalent to $\C$ and so 
$\Td\ast_\C \Td$ is $KK$-equivalent to $\C\ast_\C \C=\C$. 
\end{rem} 
 
This Lemma entails that 
the $K_1$-injectivity question for unital properly infinite \cst-algebras boils down 
to knowing whether $\U( \Td\ast_\C\Td)=\U^0( \Td\ast_\C\Td)\,$. 
Note that Proposition \ref{prop3.2} already yields that 
$\U(\Td)\ast_\T\U(\Td)\subset \U^0( \Td\ast_\C\Td)\,$. 

\bigskip 
But the following holds. 
\begin{prop}\label{prop4.5} 
{Set $p_3=\!1-s_1^{}s_1^*-s_2^{}s_2^*$  in the Cuntz algebra $\Td$
and let $u$ be the canonical unitary generating $C^*(\Z)$.} \\ 
\textrm{(i)} The relations $\jmath_0(s_k) \mapsto s_k$ and $\jmath_1(s_k) \mapsto u\,s_k$ ($k=1, 2$) 
uniquely define a unital $*$-homomorphism $\Td\ast_\C\Td\rightarrow\Td\ast_\C C^*(\Z)$ 
which is injective but not $K_1$-surjective. \\
\textrm{(ii)} The two projections $\jmath_0(p_3)$ and $\jmath_1(p_3)$ satisfy 
$\jmath_1(p_3)\,\not\sim\, \jmath_0(p_3)$ in $ \Td\ast_\C\Td$. \\ 
\textrm{(iii)} There is no $v\in\mathcal{U}( \Td\ast_\C\Td)$ such that 
\hbox{
$\jmath_1( s_1^{}s_1^*+s_2^{}s_2^*)=v\,\jmath_0(s_1^{}s_1^*+s_2^{}s_2^*)\,v^*$. 
}\\
\textrm{(iv)} There is a unitary $v\in\mathcal{U}( \Td\ast_\C\Td)$ 
such that $\jmath_1(s_1^{}s_1^*)=v\, \jmath_0(s_1^{}s_1^*)\, v^*$.  
\end{prop}
\pf 
\textrm{(i)} The unital \cst-subalgebra of $\mathcal{O}_3$ 
generated by the two isometries $s_1$ and $s_2$ is isomorphic to $\Td$, 
whence a unital \cst-embedding 
$\Td\ast_\C\Td\subset \mathcal{O}_3\ast_\C\mathcal{O}_3$ (\cite{ADEL}). 
Let $\Phi$ be the $*$-homomorphism from $\mathcal{O}_3\ast_\C \mathcal{O}_3$ to 
the  free product $\mathcal{O}_3 \ast_\C C^*(\Z) = C^*\bigl( s_1, s_2, s_3, u \bigr)$ 
fixed by the relations 
\begin{center} 
$\Phi(\jmath_0(s_k) )=s_k$ \quad and \quad $\Phi(\jmath_1(s_k) )=u\,s_k$\quad for $k=1, 2, 3$ 
\end{center}  
and let $\Psi: \mathcal{O}_3 \ast_\C C^*(\Z)\to \mathcal{O}_3\ast_\C \mathcal{O}_3$ be 
the only $*$-homomorphism such that 
\begin{center} 
$\Psi(u)=\mathop{\sum}\limits_{l=1}^3 \jmath_1(s_l^{}) \jmath_0(s_l)^*$ \quad and \quad 
$\Psi(s_k)=\jmath_0(s_k)$\quad  for $k=1, 2, 3$. 
\end{center} 

For all $k=1, 2, 3$,  we have the equalities: 

\noindent 
${}$\hspace{30pt} -- $\Psi\circ\Phi(\jmath_0(s_k) )=\Psi(s_k)=\jmath_0(s_k)\,$,
\\ ${}$\hspace{30pt} -- $\Psi\circ\Phi(\jmath_1(s_k) )=\Psi(u s_k)=\jmath_1(s_k)\,$,
\\ ${}$\hspace{30pt} -- $\Phi\circ\Psi(s_k)=\Phi(\jmath_0(s_k) )=s_k\,$. 

\noindent 
Also, $\Psi(u)^*\Psi(u)=
\sum_{l, l'}\jmath_0(s_{l'})\jmath_1(s_{l'})^*\jmath_1(s_l)\jmath_0(s_l)^*=
1_{\mathcal{O}_3\ast_\C \mathcal{O}_3}=\Psi(u)\Psi(u)^*$, 
\ie $\Psi(u)$ is a unitary in $\mathcal{O}_3\ast_\C \mathcal{O}_3$ which satisfies: 

\noindent ${}$\hspace{30pt} -- $\Phi\circ\Psi(u)=\sum_{l=1, 2, 3}\Phi(\jmath_1(s_l) )\Phi(\jmath_0(s_l)^*)=u\,$.

\noindent
Thus, $\Phi$ is an invertible unital $*$-homomorphism with inverse $\Psi$ (\cite{Blac}), 
and the restriction of $\Phi$ to the \cst-subalgebra $\Td\ast_\C \Td$ takes values 
in $\Td\ast_\C C^*(\Z)\subset\mathcal{O}_3 \ast_\C C^*(\Z)$. 

Now, there is (see \cite{G}) a six-term cyclic exact sequence 
$$\begin{array}{ccccc}
K_0(\C)=\Z&\hookrightarrow& K_0\bigl( \Td\oplus C^*(\Z)\bigr)
=\Z\oplus\Z&\rightarrow&K_0(\Td\ast_\C C^*(\Z)\,)\\
\uparrow&&&&\downarrow\\
K_1(\Td\ast_\C C^*(\Z)\,)&\leftarrow &K_1\bigl( \Td\oplus C^*(\Z)\bigr)=0\oplus \Z&\leftarrow&K_1(\C)=0
\end{array}$$ 
and so, $K_1(\Td\ast_\C C^*(\Z)\,)=\Z$, 
whereas $K_1( \Td\ast_\C\Td)=0$ by Lemma~\ref{K-Td}.

\medskip\noindent 
\textrm{(ii)} Let $\pi_0:\Td\to L(H)$ be a unital $\ast$-representation on a separable Hilbert space $H$ such that $\pi_0(p_3)$ is a rank one projection, let $\pi_1:\Td\to L(H)$ be a unital $\ast$-representation such that $\pi_1(p_3)$ is a rank two projection 
and consider the induced unital $\ast$-representation $\pi=\pi_0\ast\pi_1$ of the unital free product $\Td\ast_\C\Td$.  

Then the two projections $\pi[\jmath_0(p_3)] = \pi_0(p_3)$ and 
$\pi[\jmath_1(p_3)] = \pi_1(p_3)$ have distinct ranks and so cannot be equivalent in $L(H)$. 
Hence, $\jmath_0(p_3)\not\sim\jmath_1(p_3)$ in $\Td\ast_\C\Td$.

\medskip\noindent 
\textrm{(iii)} 
This is just a rewriting of the previous assertion since $s_1^{}s_1^*+s_2^{}s_2^*=1-p_3$. 
Indeed, the partial isometry $b=\jmath_1(s_1)\jmath_0(s_1)^*+\jmath_1(s_2)\jmath_0(s_2)^*$ defines a Murray-von Neumann equivalence in $ \Td\ast_\C\Td$ between 
the projections $\jmath_0( s_1^{}s_1^*+s_2^{}s_2^*)=1 - \jmath_0(p_3)$ and 
$\jmath_1( s_1^{}s_1^*+s_2^{}s_2^*)=1 - \jmath_1(p_3)$. 
Thus, they are unitarily equivalent in $\Td\ast_\C\Td$ if and only if 
$\jmath_0(p_3) \sim\jmath_1(p_3 )$ in $\Td\ast_\C \Td$ (\cite[Proposition~2.2.2]{LLR}).

\medskip\noindent 
\textrm{(iv)} There exists a unitary $v\in\mathcal{U}( \Td\ast_\C\Td)$ 
(which is necessarily $K_1$-trivial by Lemma~\ref{K-Td}) 
such that $\jmath_1(s_1^{}s_1^*)=v\, \jmath_0(s_1^{}s_1^*)\, v^*$.  
Indeed, we have the inequalities 
$$1 >  s_2^{}s_2^*+p_3 > s_2^{}s_2^* > 
s_2s_1(s_2^{}s_2^*+p_3)s_1^*s_2^* + s_2s_2(s_2^{}s_2^*+p_3)s_2^*s_2^*
\qquad\mathrm{in}\;\Td\,.$$ 
Thus, if we set $w:=\jmath_1(s_1)\jmath_0(s_1)^*$, then 
$1-w^*w=\jmath_0(s_2^{}s_2^*+p_3)$ and $1-ww^*=\jmath_1(s_2^{}s_2^*+p_3)$ 
are two properly infinite and full $K_0$-equivalent projections in $ \Td\ast_\C\Td$. 
Thus, there is a partial isometry $a\in\Td\ast_\C\Td$ 
with $a^*a=1-w^*w^{}$ and $aa^*=1-w^{}w^*$ (\cite{Cu81}). 
The sum $v=a+w$ has the required properties (\cite[Lemma~2.4]{BRR}). 
\qed

\begin{rems}\label{Final Rems} 
(i) The equivalence (iv)$\Leftrightarrow$(v) in Proposition~\ref{main} implies that 
all unital properly infinite \cst-algebras are $K_1$-injective if and only if 
the unitary $v\in\mathcal{U}(\Td\ast_\C\Td)$ constructed in Proposition~\ref{prop4.5}.(iv) belongs to the connected component $\mathcal{U}^0(\Td\ast_\C\Td)$. 

Note that $v \oplus 1\sim_h 1\oplus 1$ 
in $\mathcal{U}(M_2(\Td\ast_\C \Td ) )$ by \cite[Exercice 8.11]{LLR}.

\smallskip 
\noindent (ii) Let $\sigma\in\mathcal{U}(\Td)$ be the symmetry 
$\sigma=s_1^{}s_2^*+s_2^{}s_1^*+p_3\,$,  
let $v\in\mathcal{U}(\Td\ast_\C\Td)$ be a unitary such that 
$\jmath_1(s_1^{}s_1^*)=v \jmath_0(s_1^{}s_1^*) v^*$ (Proposition~\ref{prop4.5}.(iv)) and 
set $z:=v^*\jmath_1(\sigma)v\jmath_0(\sigma)\,$. 

Then, $q_1=\jmath_0(s_1^{}s_1^*)$, $q_2=\jmath_0(s_2^{}s_2^*)$ and 
$q_3=z\jmath_0(s_2^{}s_2^*)z^*$ are three properly infinite full projections in $\Td\ast_\C\Td$ which satisfy: 

\noindent
-- $q_1 q_3=\jmath_0(s_1^{}s_1^*)\, v^*\,\jmath_1(s_2^{}s_2^*)\,v=v^*\,\jmath_1(s_1^{}s_1^*) \jmath_1(s_2^{}s_2^*)\,v=0=q_1 q_2\,$, \\
\noindent 
-- $q_2 \sim_h q_1 \sim_h q_3$ in $\Td\ast_\C\Td$  
since $\sigma\in\mathcal{U}^0(\Td)$ and so $z\sim_h v^* v^{}=1$ in $\mathcal{U}(\Td\ast_\C\Td)\,$,\\
\noindent 
-- $q_1+q_3=v^*\jmath_1(s_1^{}s_1^*+s_2^{}s_2^*)v \not\sim_u 
\jmath_0(s_1^{}s_1^*+s_2^{}s_2^*)=q_1+q_2$ 
in $\Td\ast_\C\Td$ by Proposition~\ref{prop4.5}.(iii).  
\end{rems}

\medskip 
\noindent
\href{mailto:Etienne.Blanchard@math.jussieu.fr}{Etienne.Blanchard@math.jussieu.fr}\\
\address{IMJ,\quad 175, rue du Chevaleret,\quad F--75013 Paris}

\end{document}